\documentclass[a4paper,9pt]{article}
\usepackage{amsmath, amsfonts, amssymb}
\usepackage[english, russian]{babel}
\usepackage{graphicx}
\numberwithin{equation}{section}
\begin{document}
\begin{center}
\textit{Ergashev T.G.}
\end{center}

\begin{center}
THE DIRICHLET PROBLEM FOR ELLIPTIC EQUATION WITH SEVERAL SINGULAR COEFFICIENTS
\end{center}

\textbf{Abstract.} Recently found all the fundamental solutions of a
multidimensional singular elliptic equation are expressed in terms of the
well-known Lauricella hypergeometric function in many variables. In this
paper, we find a unique solution of the Dirichlet problem for an elliptic equation with several
singular coefficients in explicit form. When finding a solution, we use
decomposition formulas and some adjacent relations for the Lauricella
hypergeometric function in many variables, as well as the values of some
multidimensional improper integrals.



\section{Introduction}

It is known that the theory of boundary value problems for degenerate
equations and equations with singular coefficients is one of the rapidly
developing parts of the modern theory of partial differential equations,
which is encountered in solving many important questions of an applied
nature, for example, \cite{{A1},{A2}}. A detailed bibliography and summary of studies
of the basic boundary-value equations for degenerate equations of various
types, in particular, for elliptic equations with singular coefficients, can
be found in monographs \cite{{A4},{A3},{A5},{A6}}. In addition, generalized axisymmetric
potentials have been studied using various methods \cite{{A7},{A14},{A8},{A9},{A10},{A11},{A12},{A13}}. Omitting a huge
bibliography in which various local and non-local boundary-value problems
for mixed-type equations containing elliptic equations with singular
coefficients are studied, we note some papers which are close to the present
work. In the work \cite{A15}, fundamental solutions were constructed for the
bi-axially symmetric Helmholtz equation, and in \cite{{A33},{A16}, {A34}}  the explicit
solutions of the Dirichlet and Dirichlet-Neumann problems in one quarter of
a circle was found.

Dirichlet and Dirichlet-Neumann problems for elliptic equation with one
singular coefficient in some part of ball were investigated by Agostinelli
\cite{A25} and Olevskii \cite{A26}. Recently, Nazipov published a paper devoted to the
investigation of the Tricomi problem in a mixed domain consisting of
hemisphere and cone \cite{A17}. Fundamental solutions for the following
three-dimensional elliptic equations with two and three singular
coefficients

\begin{equation}
\label{eq1}
u_{xx} + u_{yy} + u_{zz} + {\frac{{2\alpha}} {{x}}}u_{x} + {\frac{{2\beta
}}{{y}}}u_{y} = 0,\,\,0 < 2\alpha ,\,2\beta < 1
\end{equation}

\noindent
and

\begin{equation}
\label{eq2}
u_{xx} + u_{yy} + u_{zz} + {\frac{{2\alpha}} {{x}}}u_{x} + {\frac{{2\beta
}}{{y}}}u_{y} + {\frac{{2\gamma}} {{z}}}u_{z} = 0,\,\,0 < 2\alpha ,\,2\beta
,2\gamma < 1
\end{equation}

\noindent
were constructed, respectively, in \cite{A18} and \cite{A19}. For equations (\ref{eq1}) and
(\ref{eq2}), the Dirichlet, Neumann and Holmgren problems \cite{{A20},{A21},{A22}}  were solved in
some parts of the ball.

In this paper, we study the Dirichlet problem for the equation

\begin{equation}
\label{eq2222}
H_{\alpha _{1} ,...,\alpha _{n}} ^{m,n} \left( {u} \right) \equiv
{\sum\limits_{i = 1}^{m} {u_{x_{i} x_{i}}} }   + {\sum\limits_{k = 1}^{n}
{{\frac{{2\alpha _{k}}} {{x_{k}}} }u_{x_{k}}} }   = 0,
\end{equation}

\noindent
where $m \ge 2,0 < n \le m; \quad \alpha _{k} $ are constants with $0 <
2\alpha _{k} < 1.$ Hereinafter in the present work, unless there are other
reservations, the natural number $k$ will vary from 1 to $n$, inclusive.

\section{Preliminaries}

Below we give some formulas for Euler gamma-function, Gauss hypergeometric
function, multiple Lauricella hypergeometric function (that is, Lauricella
hypergeometric function in several variables), which will be used in the
next sections.

It is known that the Euler gamma-function $\Gamma \left( {a} \right)$ has
properties \cite[pp. 17-19, (2), (10), (15)]{A29}

\begin{equation}
\label{eq3}
\Gamma \left( {a + m} \right) = \Gamma \left( {a} \right)\left( {a}
\right)_{m} ;
\quad
\Gamma \left( {a + {\frac{{1}}{{2}}}} \right) = {\frac{{\sqrt {\pi}  \Gamma
\left( {2a} \right)}}{{2^{2a - 1}\Gamma \left( {a} \right)}}},
\quad
\Gamma \left( {{\frac{{1}}{{2}}}} \right) = \sqrt {\pi}  .
\end{equation}

Here $\left( {a} \right)_{m} $ is a Pochhammer symbol, for which an equality

\begin{equation}
\label{eq4}
\left( {a} \right)_{m + n} = \left( {a} \right)_{m} \left( {a + m}
\right)_{n}
\end{equation}

\noindent
is true \cite[p.67, (5)]{A29}.

A function

\[
F\left( {a,b;c;x} \right) \equiv F{\left[ {{\begin{array}{*{20}c}
 {a,b;} \hfill \\
 {c;} \hfill \\
\end{array}} x} \right]}: = {\sum\limits_{m = 0}^{\infty}  {{\frac{{\left(
{a} \right)_{m} \left( {b} \right)_{m}}} {{\left( {c} \right)_{m}
m!}}}x^{m},\,\,\,{\left| {x} \right|} < 1\,}}
\]

\noindent
is known as the Gauss hypergeometric function and an equality

\begin{equation}
\label{eq5}
F\left( {a,b;c;1} \right) = {\frac{{\Gamma \left( {c} \right)\Gamma \left(
{c - a - b} \right)}}{{\Gamma \left( {c - a} \right)\Gamma \left( {c - b}
\right)}}}\,\,\,\,\,\,{\left[ {c \ne 0, - 1, - 2,...;\,\,Re\left( {c - a -
b} \right) > 0} \right]}
\end{equation}

\noindent
holds \cite[c.73, (14)]{A29}. Moreover, the following autotransformer formula \cite[
p.76, (22)]{A29}

\begin{equation}
\label{eq6}
F\left( {a,b;c;x} \right) = \left( {1 - x} \right)^{ - b}F\left( {c -
a,b;c;{\frac{{x}}{{x - 1}}}} \right)
\end{equation}

\noindent
is valid.

The Lauricella hypergeometric function in $n$ variables has a form \cite{A30}

\[
F_{A}^{(n)} \left( {a,b_{1} ,...,b_{n} ;c_{1} ,...,c_{n} ;x_{1} ,...,x_{n}}
\right) \equiv F_{A}^{(n)} {\left[ {{\begin{array}{*{20}c}
 {a,b_{1} ,...,b_{n} ;} \hfill \\
 {c_{1} ,...,c_{n} ;} \hfill \\
\end{array}} x_{1} ,...,x_{n}}  \right]}
\]

\begin{equation}
\label{eq7}
 = {\sum\limits_{m_{1} ,...m_{n} = 0}^{\infty}  {{\frac{{\left( {a}
\right)_{m_{1} + ... + m_{n}}  \left( {b_{1}}  \right)_{m_{1}}  ...\left(
{b_{n}}  \right)_{m_{n}}} } {{\left( {c_{1}}  \right)_{m_{1}}  ...\left(
{c_{n}}  \right)_{m_{n}}} } }{\frac{{x_{1}^{m_{1}}} } {{m_{1}
!}}}...{\frac{{x_{n}^{m_{n}}} } {{m_{n} !}}}}}
\quad
\,{\left[ {c_{k} \ne 0,\,\,\,k = \overline {1,n} ;\,\,\,{\left| {x_{1}}
\right|} + ... + {\left| {x_{n}}  \right|} < 1} \right]}.
\end{equation}

For a given multiple hypergeometric function, it is useful to fund a
decomposition formula which would express the multivariable hypergeometric
function in terms of products of several simpler hypergeometric functions
involving fewer variables. Burchnall and Chaundy \cite{{A31},{A32}} systematically
presented a number of expansion and decomposition formulas for some double
hypergeometric functions in series of simpler hypergeometric functions. For
example, the Appell function

\[
F_{2} \left( {a,b_{1} ,b_{2} ;c_{1} ,c_{2} ;x,y} \right) = {\sum\limits_{m,n
= 0}^{\infty}  {{\frac{{\left( {a} \right)_{m + n} \left( {b_{1}}
\right)_{m} \left( {b_{2}}  \right)_{n}}} {{\left( {c_{1}}  \right)_{m}
\left( {c_{2}}  \right)_{n}}} }{\frac{{x^{m}}}{{m!}}}{\frac{{y^{n}}}{{n!}}}}
}
\quad
{\left[ {c_{1} ,c_{2} \ne 0;\,\,{\left| {x} \right|} + {\left| {y} \right|}
< 1} \right]}
\]

\noindent
has the expansion \cite{A31}

\[
F_{2} \left( {a,b_{1} ,b_{2} ;c_{1} ,c_{2} ;x,y} \right) = {\sum\limits_{i =
0}^{\infty}  {{\frac{{\left( {a} \right)_{i} \left( {b_{1}}  \right)_{i}
\left( {b_{2}}  \right)_{i}}} {{i!\left( {c_{1}}  \right)_{i} \left( {c_{2}
} \right)_{i}}} }x^{i}y^{i}F\left( {a + i,b_{1} + i;c_{1} + i;x}
\right)F\left( {a + i,b_{2} + i;c_{2} + i;y} \right)}} .
\]

The Birchnell-Cendi method, which is limited to functions of two variables,
is based on the following mutually inverse symbolic operators \cite{A31}

\begin{equation}
\label{eq8}
\nabla \left( {h} \right) = {\frac{{\Gamma \left( {h} \right)\Gamma \left(
{{\rm \delta} _{1} + {\rm \delta} _{2} + h} \right)}}{{\Gamma \left( {{\rm
\delta} _{1} + h} \right)\Gamma \left( {{\rm \delta} _{2} + h} \right)}}},
\quad
\Delta \left( {h} \right) = {\frac{{\Gamma \left( {{\rm \delta} _{1} + h}
\right)\Gamma \left( {{\rm \delta} _{2} + h} \right)}}{{\Gamma \left( {h}
\right)\Gamma \left( {{\rm \delta} _{1} + {\rm \delta} _{2} + h}
\right)}}},
\end{equation}

\noindent
where ${\rm \delta} _{1} = x{\frac{{\partial}} {{\partial x}}}$ è ${\rm
\delta} _{2} = y{\frac{{\partial}} {{\partial y}}}$.

In order to generalize the operators $\nabla \left( {h} \right)$ and $\Delta
\left( {h} \right)$, defined in (\ref{eq8}), A.Hasanov and H.M.Srivastava \cite{{A23}, {A24}}
introduced the operators

\begin{equation}
\label{eq9}
\tilde {\nabla} _{z_{1} ;z_{2} ,...,z_{n}}  \left( {h} \right) =
{\frac{{{\rm \Gamma} \left( {h} \right){\rm \Gamma} \left( {{\rm \delta
}_{1} + ... + {\rm \delta} _{n} + h} \right)}}{{{\rm \Gamma} \left( {{\rm
\delta} _{1} + h} \right){\rm \Gamma} \left( {{\rm \delta} _{2} + ... + {\rm
\delta} _{n} + h} \right)}}},
\end{equation}

\begin{equation}
\label{eq10}
\tilde {\Delta} _{z_{1} ;z_{2} ,...,z_{n}}  \left( {h} \right) =
{\frac{{{\rm \Gamma} \left( {{\rm \delta} _{1} + h} \right){\rm \Gamma
}\left( {{\rm \delta} _{2} + ... + {\rm \delta} _{n} + h} \right)}}{{{\rm
\Gamma} \left( {h} \right){\rm \Gamma} \left( {{\rm \delta} _{1} + ... +
{\rm \delta} _{n} + h} \right)}}},
\end{equation}

\noindent
where ${\rm \delta} _{k} = z_{k} {\frac{{\partial}} {{\partial z_{k}}} }\,,$
with the help of which they managed to find decomposition formulas for a
whole class of hypergeometric functions in several variables. For example,
the hypergeometric Lauricelli function $F_{A}^{\left( {n} \right)} $,
defined by formula (\ref{eq7}) has the decomposition formula \cite{A23}

\[
F_{A}^{(n)} \left( {a,b_{1} ,...,b_{n} ;c_{1} ,...,c_{n} ;z_{1} ,...,z_{n}}
\right)
\]
\[
= {\sum\limits_{m_{2} ,...,m_{n} = 0}^{\infty}  {}} {\frac{{\left(
{a} \right)_{m_{2} + ... + m_{n}}  \left( {b_{1}}  \right)_{m_{2} + ... +
m_{n}}  \left( {b_{2}}  \right)_{m_{2}}  ...\left( {b_{n}}  \right)_{m_{n}}
}}{{m_{2} !...m_{n} !\left( {c_{1}}  \right)_{m_{2} + ... + m_{n}}  \left(
{c_{2}}  \right)_{m_{2}}  ...\left( {c_{n}}  \right)_{m_{n}}
}}}z_{1}^{m_{2} + ... + m_{n}}  z_{2}^{m_{2}}  ...z_{n}^{m_{n}}
\]

\[
 \cdot F\left( {a + m_{2} + ... + m_{n} ,b_{1} + m_{2} + ... + m_{n} ;c_{1}
+ m_{2} + ... + m_{n} ;z_{1}}  \right)
\]

\begin{equation}
\label{eq11}
 \cdot F_{A}^{(n - 1)} \left( {a + m_{2} + ... + m_{n} ,b_{2} + m_{2}
,...,b_{n} + m_{n} ;c_{2} + m_{2} ,...,c_{n} + m_{n} ;z_{2} ,...,z_{n}}
\right),n \in {\rm N}\backslash {\left\{ {1} \right\}}.
\end{equation}

However, due to the recurrence of formula (\ref{eq11}), additional difficulties may
arise in the applications of this expansion. Further study of the properties
of operators (\ref{eq9}) and (\ref{eq10}) showed that formula (\ref{eq11}) can be reduced to a
more convenient form.

\textbf{Lemma 1} \cite{A28}. The following decomposition formula holds true at $n
\in {\rm N}\backslash {\left\{ {1} \right\}}$

\[
F_{A}^{(n)} \left( {a,b_{1} ,b_{2} ,....,b_{n} ;c_{1} ,c_{2} ,....,c_{n}
;z_{1} ,...,z_{n}}  \right)
\]

\begin{equation}
\label{eq12}
 = {\sum\limits_{{\mathop {m_{i,j} = 0}\limits_{(2 \le i \le j \le n)}
}}^{\infty}  {{\frac{{(a)_{N(n,n)}}} {{m_{ij} !}}}}} {\prod\limits_{k =
1}^{n} {{\left[ {{\frac{{(b_{k} )_{M(k,n)}}} {{(c_{k} )_{M(k,n)}
}}}z_{k}^{M(k,n)} F\left( {a + N(k,n),b_{k} + M(k,n);c_{k} + M(k,n);z_{k}}
\right)} \right]}}} ,
\end{equation}

\noindent
where

\begin{equation}
\label{eq13}
M(k,n) = {\sum\limits_{i = 2}^{k} {m_{i,k} +}}  {\sum\limits_{i = k + 1}^{n}
{m_{k + 1,i}}}  ,
\quad
N(k,n) = {\sum\limits_{i = 2}^{k + 1} {{\sum\limits_{j = i}^{n} {m_{i,j}}
}}} .
\end{equation}

The formula (\ref{eq12}) is proved by the method mathematical induction \cite{A28}.

It should be noted here that the sum ${\sum\limits_{k = 1}^{n} {M(k,n)}} $
has the parity property, which plays an important role in the calculation of
the some values of hypergeometric functions. In fact, by virtue of equality

\[
{\sum\limits_{k = 2}^{n} {{\sum\limits_{i = 2}^{k} {m_{i,k}}} } }  =
{\sum\limits_{k = 1}^{n - 1} {{\sum\limits_{i = k + 1}^{n} {m_{k + 1,i}}} }
}
\]

\noindent
we obtain

\begin{equation}
\label{eq14}
{\sum\limits_{k = 1}^{n} {M(k,n)}}  = 2{\sum\limits_{k = 2}^{n}
{{\sum\limits_{i = 2}^{k} {m_{i,k}}} } }  = 2{\sum\limits_{k = 1}^{n - 1}
{{\sum\limits_{i = k + 1}^{n} {m_{k + 1,i}}} } } .
\end{equation}

In the present paper, ${\rm R}_{m}^{n +}  $ denotes $1 / 2^{n}$ part of the
Euclidean space ${\rm R}^{m}$:

\[
{\rm R}_{m}^{n +}  : = {\left\{ {\left( {x_{1} ,...,x_{m}}
\right):\,\,x_{1} > 0,\,...,\,x_{n} > 0,\,\,1 \le n \le m,\,\,m \ge 2}
\right\}}.
\]

All the fundamental solutions of equation (\ref{eq2222}) in the domain ${\rm
R}_{m}^{n +}  $ were found in \cite{A28}, and we will use one of these solutions
in the study of the problem:

\begin{equation}
\label{eq15}
q_{n} \left( {{\rm x};{\rm \xi}}  \right) = \gamma _{n}
{\prod\limits_{i = 1}^{n} {{\left[ {\left( {x_{i} \xi _{i}}  \right)^{1 -
2\alpha _{i}}}  \right]} \cdot}}  \,r^{ - 2\tilde {\alpha} _{n}
}F_{A}^{\left( {n} \right)} \left( {\tilde {\alpha} _{n} ,1 - \alpha _{1}
,...,1 - \alpha _{n} ;2 - 2\alpha _{1} ,...,2 - 2\alpha _{n} ;\sigma}
\right),
\end{equation}

\noindent
where

\[
{\rm x}: = \left( {x_{1} ,...,x_{m}}  \right),
{\rm \xi} : = \left( {\xi _{1} ,...,\xi _{m}}  \right),
\quad
{\rm \sigma} : = \left( {\sigma _{1} ,...,\sigma _{n}}  \right);
\]

\[
0 < 2\alpha _{1} ,...,2\alpha _{n} < 1,
\quad
\tilde {\alpha} _{n} = {\frac{{m - 2}}{{2}}} + n - \alpha _{1} - ... -
\alpha _{n} ;
\]

\begin{equation}
\label{eq16}
\gamma _{n} = 2^{2\tilde {\alpha} _{n} - m}{\frac{{\Gamma \left( {\tilde
{\alpha} _{n}}  \right)}}{{\pi ^{m / 2}}}}{\prod\limits_{j = 1}^{n}
{{\frac{{\Gamma \left( {1 - \alpha _{j}}  \right)}}{{\Gamma \left( {2 -
2\alpha _{j}}  \right)}}}}} ,
\end{equation}

\[
r^{2} = {\sum\limits_{i = 1}^{m} {\left( {x_{i} - \xi _{i}}  \right)^{2}}} ,
\quad
r_{k}^{2} = \left( {x_{k} + \xi _{k}}  \right)^{2} + {\sum\limits_{i = 1,i
\ne k}^{m} {\left( {x_{i} - \xi _{i}}  \right)^{2}}} ,
\quad
\sigma _{k} = 1 - {\frac{{r_{k}^{2}}} {{r^{2}}}}.
\]

It is easy to verify that the fundamental solution $q_{n} \left( {{\rm
x};{\rm \xi}}  \right)$ has the property

\[
{\left. {q_{n} \left( {{\rm x};{\rm \xi}}  \right)} \right|}_{x_{k} =
0} = 0.
\]

\section{Formulation of the problem and the uniqueness of the solution}

Let $\Omega \subset {\rm R}_{m}^{n +}  $ be a finite simple-connected domain
bounded by planes $x_{1} = 0,...,x_{n} = 0$ and by smooth $m -
$dimensional surface $S$. The intersection of this surface with plane $x_{k}
= 0$ is denoted by $\chi _{k} $. Designate as the domain $S_{k} ,$ a
hyperplane $Ox_{1} ...x_{k - 1} x_{k + 1} ...x_{m} ,$ bounded by $x_{k} = 0$
($0 < x_{l} < a_{l} , \quad  - b_{s} < x_{s} < c_{s} , \quad l = \overline {1,n}
,\,\,l \ne k;\,\,n < s \le m)$ and by a curve $\chi _{k} $. Here $a_{l}
,\,\,b_{s} $ and $c_{s} $ are positive constants. We introduce the notation:

\[
\,\tilde {x}_{k} : = \left( {x_{1} ,...,x_{k - 1} ,x_{k + 1} ,...,x_{n}
,...,x_{m}}  \right) \in S_{k} \subset {\rm R}_{m - 1}^{(n - 1) +}  \subset
{\rm R}^{m - 1}.
\]

\textbf{Dirichlet problem.} To find a function $u\left( {{\rm x}} \right)
\in C\left( {\bar {\Omega}}  \right) \cap C^{2}\left( {\Omega}  \right)$,
satisfying equation (\ref{eq2222}) in $\Omega $ and conditions

\begin{equation}
\label{eq17}
{\left. {u} \right|}_{x_{k} = 0} = \tau _{k} \left( {\tilde {x}_{k}}
\right),
\quad
\,\tilde {x}_{k} \in \bar {S}_{k} ,
\end{equation}

\begin{equation}
\label{eq18}
{\left. {u} \right|}_{S} = \varphi \left( {{\rm x}} \right),
\quad
\,{\rm x} \in \bar {S},
\end{equation}

\noindent
where $\tau _{k} \left( {\tilde {x}_{k}}  \right)$ and $\varphi \left( {{\rm
x}} \right)$ are given continuous functions fulfilling the following
matching conditions:

\begin{equation}
\label{eq19}
\tau _{i} \left( {0,...,0,x_{n + 1} ,...,x_{m}}  \right) = \tau _{j} \left(
{0,...,0,x_{n + 1} ,...,x_{m}}  \right),
i \ne j,
i,\,j = \overline {1,n} ,
\quad
{\left. {\tau _{k} \left( {\tilde {x}_{k}}  \right)} \right|}_{\chi _{k}}  =
{\left. {\varphi \left( {{\rm x}} \right)} \right|}_{\chi _{k}}  .
\end{equation}

One can readily check the validity of the following relation

\[
x_{1}^{2\alpha _{1}}  ...x_{n}^{2\alpha _{n}}  {\left[ {uH_{\alpha} ^{m,n}
\left( {w} \right) - wH_{\alpha} ^{m,n} \left( {u} \right)} \right]} =
x_{1}^{2\alpha _{1}}  ...x_{n}^{2\alpha _{n}}  {\sum\limits_{i = 1}^{m}
{{\frac{{\partial}} {{\partial x_{i}}} }}} \left( {uw_{x_{i}}  - wu_{x_{i}}
} \right),\,\,n \le m.
\]

Let $\Omega _{\varepsilon}  $ be a sub-domain of $\Omega $ at a distance
$\varepsilon > 0$ from its boundary $\partial \Omega  =
{\bigcup\limits_{i = 1}^{n} {S_{i}}}   \cup S$ and

\[
{\frac{{\partial}} {{\partial {\rm {\bf n}}}}} = {\sum\limits_{i = 1}^{m}
{\cos \left( {{\rm {\bf n}},x_{i}}  \right)}}  \cdot {\frac{{\partial
}}{{\partial x_{i}}} },
\]
${\rm {\bf n}}$ is outer normal to $\partial \Omega $.

Integrate both sides of above given equality on the domain $\Omega
_{\varepsilon}  $ and use the classical formula of Gauss-Ostrogadsky:

\[
{\int_{\Omega _{\varepsilon}}   {x_{1}^{2\alpha _{1}}  ...x_{n}^{2\alpha
_{n}}} }  {\left[ {uH_{\alpha _{1} ,...,\alpha _{n}} ^{m,n} \left( {w}
\right) - wH_{\alpha _{1} ,...,\alpha _{n}} ^{m,n} \left( {u} \right)}
\right]}dx_{1} ...dx_{m}
\]

\begin{equation}
\label{eq20}
 = {\int_{\partial \Omega _{\varepsilon}}   {x_{1}^{2\alpha _{1}}
...x_{n}^{2\alpha _{n}}  {\sum\limits_{i = 1}^{m} {\left( {uw_{x_{i}}  -
wu_{x_{i}}}   \right)\cos \left( {{\rm {\bf n}},x_{i}}  \right)d\vartheta}
}}} .
\end{equation}

Using the equality

\[
x_{1}^{2\alpha _{1}}  ...x_{n}^{2\alpha _{n}}  {\left[ {uH_{\alpha _{1}
,...,\alpha _{n}} ^{m,n} \left( {u} \right) + {\sum\limits_{i = 1}^{m}
{\left( {{\frac{{\partial u}}{{\partial x_{i}}} }} \right)^{2}}}}  \right]}
= {\sum\limits_{i = 1}^{m} {{\frac{{\partial}} {{\partial x_{i}}} }}} \left(
{x_{1}^{2\alpha _{1}}  ...x_{n}^{2\alpha _{n}}  u{\frac{{\partial
u}}{{\partial x_{i}}} }} \right),
\]

\noindent
we obtain

\[
{\int_{\Omega _{\varepsilon}}   {x_{1}^{2\alpha _{1}}  ...x_{n}^{2\alpha
_{n}}} }  uH_{\alpha _{1} ,...,\alpha _{n}} ^{m,n} \left( {u} \right)dx_{1}
...dx_{m} + {\int_{\Omega _{\varepsilon}}   {x_{1}^{2\alpha _{1}}
...x_{n}^{2\alpha _{n}}} }  {\sum\limits_{i = 1}^{m} {\left(
{{\frac{{\partial u}}{{\partial x_{i}}} }} \right)^{2}}} dx_{1} ...dx_{m}
\]

\[
 = {\int_{\Omega _{\varepsilon}}   {{\sum\limits_{i = 1}^{m}
{{\frac{{\partial}} {{\partial x_{i}}} }}} \left( {x_{1}^{2\alpha _{1}}
...x_{n}^{2\alpha _{n}}  u{\frac{{\partial u}}{{\partial x_{i}}} }}
\right)dx_{1} ...dx_{m}}}  .
\]

Applying again the formula of Gauss-Ostrogradsky to this equality and
letting $\varepsilon \to 0$, we get

\[
{\int_{\Omega _{\varepsilon}}   {x_{1}^{2\alpha _{1}}  ...x_{n}^{2\alpha
_{n}}} }  {\sum\limits_{i = 1}^{m} {\left( {{\frac{{\partial u}}{{\partial
x_{i}}} }} \right)^{2}}} dx_{1} ...dx_{m}
\]

\begin{equation}
\label{eq21}
 = {\sum\limits_{k = 1}^{n} {{\int_{S_{k}}  {x_{1}^{2\alpha _{1}}  ...x_{k -
1}^{2\alpha _{k - 1}}  x_{k + 1}^{2\alpha _{k + 1}}  ...x_{n}^{2\alpha _{n}
} \tau _{k} \nu _{k} dS_{k}}} } }  + {\int_{S} {x_{1}^{2\alpha _{1}}
...x_{n}^{2\alpha _{n}}  \varphi {\frac{{\partial u}}{{\partial {\rm {\bf
n}}}}}dS,}}
\end{equation}

\noindent
where

\[
\nu _{k} \left( {\tilde {x}_{k}}  \right): = {\left. {\left( {x_{k}^{2\alpha
_{k}}  {\frac{{\partial u}}{{\partial x_{k}}} }} \right)} \right|}_{x_{k} =
0} .
\]

To prove the uniqueness of the solution, as usual, we suppose that the
problem has two $v, \quad w$ solutions. Denoting $u = v - w$ we have that
satisfies homogeneous Dirichlet problem ($\tau _{k} = 0, \quad \varphi = 0)$.
Further we have to prove that the homogeneous problem has only trivial
solution. In this case from (\ref{eq21}) one can easily get

\[
{\int_{\Omega _{\varepsilon}}   {x_{1}^{2\alpha _{1}}  ...x_{n}^{2\alpha
_{n}}} }  {\sum\limits_{i = 1}^{m} {\left( {{\frac{{\partial u}}{{\partial
x_{i}}} }} \right)^{2}}} dx_{1} ...dx_{m} = 0.
\]

Hence, it follows that $u_{x_{1}}  = ... = u_{x_{m}}  = 0,$ which implies
that $u$ is a constant function. Considering homogeneous conditions (\ref{eq17})
and (\ref{eq18}), we conclude that $u\left( {{\rm x}} \right) \equiv 0$ in
$\overline {\Omega}  $.

\section{The existence of the solution}

\bigskip

We prove the existence of the solution in a special case of the domain
$\Omega $ in order to get the solution in an explicit form. Assume $R =
a_{k} = b_{k} = c_{k} $ and let

\[
\Omega = {\left\{ {{\rm x}: x_{1}^{2} + ... + x_{m}^{2} <
R^{2},\,\,x_{1} > 0,...,x_{n} > 0,\,\, - R < x_{n + 1} < R,..., - R < x_{m}
< R} \right\}}.
\]

We find a solution of considered problem using method Green's functions \cite{A27}
. Therefore, first we give a definition of Green's function for the
formulated problem.

\textbf{Definition.} We call the function $G\left( {{\rm x};{\rm \xi}}
\right)$ as Green's function of the Dirichlet problem, if it satisfies the
following conditions:

\noindent
this function is a regular solution of equation (\ref{eq2222}) in the domain $\Omega
$, expect at the point ${\rm \xi} $, which is any fixed point of $\Omega $;

\noindent
it satisfies boundary conditions

\[
{\left. {G\left( {{\rm x};{\rm \xi}}  \right)} \right|}_{x_{k} = 0} = 0,
\quad
{\left. {G\left( {{\rm x};{\rm \xi}}  \right)} \right|}_{S} = 0;
\]

\noindent
it can be represented as

\begin{equation}
\label{eq22}
G\left( {{\rm x};{\rm \xi}}  \right) = q_{n} \left( {{\rm x};{\rm \xi}}
\right) + q_{n}^{ *}  \left( {{\rm x};{\rm \xi}}  \right),
\end{equation}

\noindent
where $q_{n} \left( {{\rm x};{\rm \xi}}  \right)$ is the fundamental
solution found earlier (see a formula (\ref{eq15})), function

\[
q_{n}^{ *}  \left( {{\rm x};{\rm \xi}}  \right) = - \left(
{{\frac{{a}}{{R_{0}}} }} \right)^{2\tilde {\alpha} _{n}} q_{n} \left( {{\rm
x};{\rm \bar {\xi}} } \right)
\]

\noindent
is a regular solution of equation (\ref{eq2222}) in the domain $\Omega $. Here

\[
{\rm \bar {\xi}} : = \left( {\bar {\xi} _{1} ,...,\bar {\xi} _{m}}  \right),
\quad
\bar {\xi} _{i} = {\frac{{a^{2}}}{{R_{0}^{2}}} }\xi _{i} ,
\quad
R_{0}^{2} = \xi _{1}^{2} + ... + \xi _{m}^{2} .
\]

Excise a small ball with its center at ${\rm \xi} $ and with radius $\rho >
0$ from the domain $\Omega $. Designate the sphere of the excised ball as
$C_{\rho}  $ and by $\Omega _{\rho}  $ denote the remaining part of $\Omega
$.

Applying formula (\ref{eq20}), we obtain

\[
{\int_{C_{\rho}}   {{\rm x}^{\left( {2\alpha}  \right)}{\left[ {u\left(
{{\rm x}} \right){\frac{{\partial G\left( {{\rm x};{\rm \xi}}
\right)}}{{\partial {\rm {\bf n}}}}} - G\left( {{\rm x};{\rm \xi}}
\right){\frac{{\partial u\left( {{\rm x}} \right)}}{{\partial {\rm {\bf
n}}}}}} \right]}dC_{\rho}} }
\]

\begin{equation}
\label{eq23}
 = {\sum\limits_{k = 1}^{n} {{\int_{S_{k}}  {\tilde {x}_{k}^{\left( {2\alpha
} \right)} G_{k}^{\ast}  \left( {x_{1} ,...,x_{k - 1} ,0,x_{k + 1}
,...,x_{m} ;{\rm \xi}}  \right)\tau _{k} \left( {\tilde {x}_{k}}
\right)dS_{k}}}  \,}}  + {\int_{S} {{\rm x}^{\left( {2\alpha}
\right)}{\frac{{\partial G\left( {{\rm x};{\rm \xi}}
\right)}}{{\partial {\rm {\bf n}}}}}\varphi \left( {\vartheta}
\right)d\vartheta}}
\end{equation}

\noindent
where

\[
{\rm x}^{\left( {2\alpha}  \right)}: = x_{1}^{2\alpha _{1}}
...x_{n}^{2\alpha _{n}}  ,
\quad
\tilde {x}_{k}^{\left( {2\alpha}  \right)} : = x_{1}^{2\alpha _{1}}  ...x_{k
- 1}^{2\alpha _{k - 1}}  x_{k + 1}^{2\alpha _{k + 1}}  ...x_{n}^{2\alpha
_{n}}  .
\]

\[
{\left. {G_{k}^{\ast}  \left( {x_{1} ,...,x_{k - 1} ,0,x_{k + 1} ,...,x_{m}
;{\rm \xi}}  \right) = \left( {x_{k}^{2\alpha _{k}}  {\frac{{\partial
G\left( {{\rm x};{\rm \xi}}  \right)}}{{\partial x_{k}}} }} \right)}
\right|}_{x_{k} = 0} ,\,\,\,\,\,\left( {x_{1} ,...,x_{k - 1} ,0,x_{k + 1}
,...,x_{m}}  \right) \in S_{k} .
\]

First, we consider an integral

\[
{\int_{C_{\rho}}   {{\rm x}^{\left( {2\alpha}  \right)}u\left( {{\rm x}}
\right){\frac{{\partial G\left( {{\rm x};{\rm \xi}}  \right)}}{{\partial
{\rm {\bf n}}}}}dC_{\rho}} }  .
\]

Taking (\ref{eq22}) into account we rewrite it as follows

\[
{\int_{C_{\rho}}   {{\rm x}^{\left( {2\alpha}  \right)}u{\frac{{\partial
G}}{{\partial {\rm {\bf n}}}}}dC_{\rho}} }   \equiv {\int_{C_{\rho}}   {{\rm
x}^{\left( {2\alpha}  \right)}u{\frac{{\partial q_{n} \left( {{\rm x};{\rm
\xi}}  \right)}}{{\partial {\rm {\bf n}}}}}dC_{\rho}} }   + {\int_{C_{\rho}
} {{\rm x}^{\left( {2\alpha}  \right)}u{\frac{{\partial q_{n}^{\ast}  \left(
{{\rm x};{\rm \xi}}  \right)}}{{\partial {\rm {\bf n}}}}}dC_{\rho}} }
= I_{1} + I_{2} .
\]

Using the formula of differentiation

\begin{equation}
\label{eq24}
\begin{array}{l}
 {\frac{{\partial}} {{\partial z_{j}}} }F_{A}^{(n)} \left( {a,b_{1}
,...,b_{n} ;c_{1} ,...,c_{n} ;z_{1} ,...,z_{n}}  \right) \\
\\
 = {\frac{{ab_{j}}} {{c_{j}}} }F_{A}^{(n)} \left( {a + 1,b_{1} ,...,b_{j -
1} ,b_{j} + 1,b_{j + 1} ,...,b_{n} ;c_{1} ,...,c_{j - 1} ,c_{j} + 1,c_{j +
1} ,...,c_{n} ;z_{1} ,...,z_{n}}  \right) \\
 \end{array}
\end{equation}

\noindent
and the following adjacent relation

\begin{equation}
\label{eq25}
\begin{array}{l}
 {\sum\limits_{j = 1}^{n} {{\frac{{a_{j}}} {{b_{j}}} }x_{j} F_{A}^{(n)}
\left( {a + 1,a_{1} ,...,a_{j - 1} ,a_{j} + 1,a_{j + 1} ,...,a_{n} ;b_{1}
,...,b_{j - 1} ,b_{j} + 1,b_{j + 1} ,...,b_{n} ;x_{1} ,...,x_{n}}  \right)}
} \\
\\
 = F_{A}^{(n)} \left( {a + 1,a_{1} ,...,a_{n} ;b_{1} ,...,b_{n} ;x_{1}
,...,x_{n}}  \right) - F_{A}^{(n)} \left( {a,a_{1} ,...,a_{n} ;b_{1}
,...,b_{n} ;x_{1} ,...,x_{n}}  \right). \\
 \end{array}
\end{equation}

\noindent
we calculate

\begin{equation}
\label{eq26}
{\frac{{\partial q_{n} \left( {{\rm x};{\rm \xi}}  \right)}}{{\partial {\rm
{\bf n}}}}} = {\sum\limits_{i = 1}^{m} {{\frac{{\partial q_{n} \left( {{\rm
x};{\rm \xi}}  \right)}}{{\partial x_{i}}} }}}  \cdot \cos \left( {{\rm {\bf
n}},\,x_{i}}  \right).
\end{equation}

Below we get detailed evaluations for ${\frac{{\partial q_{n} \left( {{\rm
x};{\rm \xi}}  \right)}}{{\partial x_{i}}} }$, when $1 \le i \le n$. Indeed,
using the formula of differentiation (\ref{eq24}), we get

\[
{\frac{{\partial q_{n} \left( {x,\xi}  \right)}}{{\partial x_{k}}} } =
\left( {1 - 2\alpha _{k}}  \right)\gamma _{n} x_{k}^{ - 2\alpha _{k}}  \xi
_{k}^{1 - 2\alpha _{k}}  {\prod\limits_{i = 1,i \ne k}^{n} {{\left[ {\left(
{x_{i} \xi _{i}}  \right)^{1 - 2\alpha _{i}}}  \right]}\,}} r^{ - 2\tilde
{\alpha} _{n}} F_{A}^{\left( {n} \right)} {\left[ {{\begin{array}{*{20}c}
 {\tilde {\alpha} _{n} ,1 - \alpha _{1} ,...,1 - \alpha _{n} ;} \hfill \\
 {2 - 2\alpha _{1} ,...,2 - 2\alpha _{n} ;} \hfill \\
\end{array}} \sigma}  \right]}
\]

\[
 - 2\tilde {\alpha} _{n} \gamma _{n} \left( {x_{k} - \xi _{k}}
\right){\prod\limits_{i = 1}^{n} {{\left[ {\left( {x_{i} \xi _{i}}
\right)^{1 - 2\alpha _{i}}}  \right]} \cdot}}  r^{ - 2\tilde {\alpha} _{n} -
2}F_{A}^{\left( {n} \right)} {\left[ {{\begin{array}{*{20}c}
 {\tilde {\alpha} _{n} ,1 - \alpha _{1} ,...,1 - \alpha _{n} ;} \hfill \\
 {2 - 2\alpha _{1} ,...,2 - 2\alpha _{n} ;} \hfill \\
\end{array}} \sigma}  \right]}
\]

\[
\begin{array}{l}
 - 2\tilde {\alpha} _{n} \gamma _{n} \xi _{k} {\prod\limits_{i = 1}^{n}
{{\left[ {\left( {x_{i} \xi _{i}}  \right)^{1 - 2\alpha _{i}}}  \right]}}}
\cdot r^{ - 2\tilde {\alpha} _{n} - 2}
\\
\\
\times F_{A}^{\left( {n} \right)} {\left[
{{\begin{array}{*{20}c}
 {\tilde {\alpha} _{n} + 1,1 - \alpha _{1} ,...,1 - \alpha _{k - 1} ,2 -
\alpha _{k} ,1 - \alpha _{k + 1} ,...1 - \alpha _{n} ;} \hfill \\
 {2 - 2\alpha _{1} ,...,2 - 2\alpha _{k - 1} ,3 - 2\alpha _{k} ,2 - 2\alpha
_{k + 1} ,...2 - 2\alpha _{n} ;} \hfill \\
\end{array}} \sigma}  \right]}
\end{array}
\]

\[
\begin{array}{l}
 - 2\tilde {\alpha} _{n} \gamma _{n} \left( {x_{k} - \xi _{k}}
\right){\prod\limits_{i = 1}^{n} {{\left[ {\left( {x_{i} \xi _{i}}
\right)^{1 - 2\alpha _{i}}}  \right]}}} \,r^{ - 2\tilde {\alpha} _{n} - 2}
\\
\\
 \times {\sum\limits_{i = 1}^{n} {}} {\frac{{1 - \alpha _{i}}} {{2 - 2\alpha
_{i}}} }\sigma _{i} F_{A}^{\left( {n} \right)} {\left[
{{\begin{array}{*{20}c}
 {\tilde {\alpha} _{n} + 1,1 - \alpha _{1} ,...,1 - \alpha _{i - 1} ,2 -
\alpha _{i} ,1 - \alpha _{i + 1} ,...1 - \alpha _{n} ;} \hfill \\
 {2 - 2\alpha _{1} ,...,2 - 2\alpha _{i - 1} ,3 - 2\alpha _{i} ,2 - 2\alpha
_{i + 1} ,...2 - 2\alpha _{n} ;} \hfill \\
\end{array}} \sigma}  \right]}. \\
 \end{array}
\]

Considering adjacent relation (\ref{eq25}) we obtain

\[
{\frac{{\partial q_{n} \left( {x,\xi}  \right)}}{{\partial x_{k}}} } =
\left( {1 - 2\alpha _{k}}  \right)\gamma _{n} x_{k}^{ - 2\alpha _{k}}  \xi
_{k}^{1 - 2\alpha _{k}}  {\prod\limits_{i = 1,i \ne k}^{n} {{\left[ {\left(
{x_{i} \xi _{i}}  \right)^{1 - 2\alpha _{i}}}  \right]}\,}} r^{ - 2\tilde
{\alpha} _{n}} F_{A}^{\left( {n} \right)} {\left[ {{\begin{array}{*{20}c}
 {\tilde {\alpha} _{n} ,1 - \alpha _{1} ,...,1 - \alpha _{n} ;} \hfill \\
 {2 - 2\alpha _{1} ,...,2 - 2\alpha _{n} ;} \hfill \\
\end{array}} \sigma}  \right]}
\]

\[
 + 2\tilde {\alpha} _{n} \gamma _{n} \left( {\xi _{k} - x_{k}}
\right){\prod\limits_{i = 1}^{n} {{\left[ {\left( {x_{i} \xi _{i}}
\right)^{1 - 2\alpha _{i}}}  \right]}}} \,r^{ - 2\tilde {\alpha} _{n} -
2}F_{A}^{\left( {n} \right)} {\left[ {{\begin{array}{*{20}c}
 {\tilde {\alpha} _{n} + 1,1 - \alpha _{1} ,...,1 - \alpha _{n} ;} \hfill \\
 {2 - 2\alpha _{1} ,...,2 - 2\alpha _{n} ;} \hfill \\
\end{array}} \sigma}  \right]}
\]

\begin{equation}
\label{eq27}
 - 2\tilde {\alpha} _{n} \gamma _{n} \xi _{k} {\prod\limits_{i = 1}^{n}
{{\left[ {\left( {x_{i} \xi _{i}}  \right)^{1 - 2\alpha _{i}}}  \right]}}
}\,r^{ - 2\tilde {\alpha} _{n} - 2}F_{A}^{\left( {n} \right)} {\left[
{{\begin{array}{*{20}c}
 {\tilde {\alpha} _{n} + 1,1 - \alpha _{1} ,...,1 - \alpha _{k - 1} ,2 -
\alpha _{k} ,1 - \alpha _{k + 1} ,...1 - \alpha _{n} ;} \hfill \\
 {2 - 2\alpha _{1} ,...,2 - 2\alpha _{k - 1} ,3 - 2\alpha _{k} ,2 - 2\alpha
_{k + 1} ,...2 - 2\alpha _{n} ;} \hfill \\
\end{array}} \sigma}  \right]}
\end{equation}

Similarly we calculate ${\frac{{\partial q_{n} \left( {{\rm x};{\rm \xi}}
\right)}}{{\partial x_{i}}} }$, when $n + 1 \le i \le m:$

\begin{equation}
\label{eq28}
{\frac{{\partial q_{n} \left( {x,\xi}  \right)}}{{\partial x_{k}}} } =
2\tilde {\alpha} _{n} \gamma _{n} \left( {\xi _{k} - x_{k}}
\right){\prod\limits_{i = 1}^{n} {{\left[ {\left( {x_{i} \xi _{i}}
\right)^{1 - 2\alpha _{i}}}  \right]}}} \,r^{ - 2\tilde {\alpha} _{n} -
2}F_{A}^{\left( {n} \right)} {\left[ {{\begin{array}{*{20}c}
 {\tilde {\alpha} _{n} + 1,1 - \alpha _{1} ,...,1 - \alpha _{n} ;} \hfill \\
 {2 - 2\alpha _{1} ,...,2 - 2\alpha _{n} ;} \hfill \\
\end{array}} \sigma}  \right]}.
\end{equation}

Taking (\ref{eq26}), (\ref{eq27}) and (\ref{eq28}) into account we calculate

\[
\begin{array}{l}
 {\frac{{\partial q_{n} \left( {x,\xi}  \right)}}{{\partial {\rm {\bf n}}}}}
= - \tilde {\alpha} _{n} \gamma _{n} {\prod\limits_{i = 1}^{n} {{\left[
{\left( {x_{i} \xi _{i}}  \right)^{1 - 2\alpha _{i}}}  \right]}}} \,r^{ -
2\tilde {\alpha} _{n}} F_{A}^{\left( {n} \right)} {\left[
{{\begin{array}{*{20}c}
 {\tilde {\alpha} _{n} + 1,1 - \alpha _{1} ,...,1 - \alpha _{n} ;} \hfill \\
 {2 - 2\alpha _{1} ,...,2 - 2\alpha _{n} ;} \hfill \\
\end{array}} \sigma}  \right]}{\frac{{\partial}} {{\partial {\rm {\bf
n}}}}}{\left[ {\ln r^{2}} \right]} \\
\\
 - 2\tilde {\alpha} _{n} \gamma _{n} {\prod\limits_{i = 1}^{n} {{\left[
{\left( {x_{i} \xi _{i}}  \right)^{1 - 2\alpha _{i}}}  \right]}}} \,r^{ -
2\tilde {\alpha} _{n} - 2}\\
\\
\times{\sum\limits_{i = 1}^{n} {\xi _{i} F_{A}^{\left(
{n} \right)} {\left[ {{\begin{array}{*{20}c}
 {\tilde {\alpha} _{n} + 1,1 - \alpha _{1} ,...,1 - \alpha _{i - 1} ,2 -
\alpha _{i} ,1 - \alpha _{i + 1} ,...1 - \alpha _{n} ;} \hfill \\
 {2 - 2\alpha _{1} ,...,2 - 2\alpha _{i - 1} ,3 - 2\alpha _{i} ,2 - 2\alpha
_{i + 1} ,...2 - 2\alpha _{n} ;} \hfill \\
\end{array}} \sigma}  \right]}\cos \left( {{\rm {\bf n}};x_{i}}  \right)}}
\\
\\
 + \gamma _{n} {\prod\limits_{i = 1}^{n} {{\left[ {\left( {x_{i} \xi _{i}}
\right)^{1 - 2\alpha _{i}}}  \right]}}} \,r^{ - 2\tilde {\alpha} _{n}
}F_{A}^{\left( {n} \right)} {\left[ {{\begin{array}{*{20}c}
 {\tilde {\alpha} _{n} ,1 - \alpha _{1} ,...,1 - \alpha _{n} ;} \hfill \\
 {2 - 2\alpha _{1} ,...,2 - 2\alpha _{n} ;} \hfill \\
\end{array}} \sigma}  \right]} \cdot {\sum\limits_{i = 1}^{n} {{\frac{{1 -
2\alpha _{i}}} {{x_{i}}} }\cos \left( {{\rm {\bf n}};x_{i}}  \right)}}  \\
 \end{array}
\]

Now consider the integral

\[
{\int_{C_{\rho}}   {{\rm x}^{\left( {2\alpha}  \right)}u{\frac{{\partial
q_{n} \left( {{\rm x};{\rm \xi}}  \right)}}{{\partial {\rm {\bf
n}}}}}dC_{\rho}} }   = I_{11} + I_{12} + I_{13} ,
\]

\noindent
where

\[
I_{11} = - \tilde {\alpha} _{n} \gamma _{n} {\prod\limits_{i = 1}^{n}
{{\left[ {\xi _{i}^{1 - 2\alpha _{i}}}   \right]}}} {\int_{C_{\rho}}   {u
\cdot {\prod\limits_{i = 1}^{n} {{\left[ {x_{i}}  \right]}}} \,r^{ - 2\tilde
{\alpha} _{n}} F_{A}^{\left( {n} \right)} {\left[ {{\begin{array}{*{20}c}
 {\tilde {\alpha} _{n} + 1,1 - \alpha _{1} ,...,1 - \alpha _{n} ;} \hfill \\
 {2 - 2\alpha _{1} ,...,2 - 2\alpha _{n} ;} \hfill \\
\end{array}} \sigma}  \right]}{\frac{{\partial}} {{\partial {\rm {\bf
n}}}}}{\left[ {\ln r^{2}} \right]}dC_{\rho}} }  ,
\]

\[
I_{12} = - 2\tilde {\alpha} _{n} \gamma _{n} {\prod\limits_{i = 1}^{n}
{{\left[ {\xi _{i}^{1 - 2\alpha _{i}}}   \right]}}} \,{\int_{C_{\rho}}   {u
\cdot {\prod\limits_{i = 1}^{n} {{\left[ {x_{i}}  \right]}\,}} r^{ - 2\tilde
{\alpha} _{n} - 2}}}
\]

\[
\times {\sum\limits_{i = 1}^{n} {\xi _{i} F_{A}^{\left( {n} \right)} {\left[
{{\begin{array}{*{20}c}
 {\tilde {\alpha} _{n} + 1,1 - \alpha _{1} ,...,1 - \alpha _{i - 1} ,2 -
\alpha _{i} ,1 - \alpha _{i + 1} ,...1 - \alpha _{n} ;} \hfill \\
 {2 - 2\alpha _{1} ,...,2 - 2\alpha _{i - 1} ,3 - 2\alpha _{i} ,2 - 2\alpha
_{i + 1} ,...2 - 2\alpha _{n} ;} \hfill \\
\end{array}} \sigma}  \right]}\cos \left( {{\rm {\bf n}};x_{i}}  \right)}
}dC_{\rho}  ,
\]

\[
I_{13} = \gamma _{n} {\prod\limits_{i = 1}^{n} {{\left[ {\xi _{i}^{1 -
2\alpha _{i}}}   \right]}}} {\int_{C_{\rho}}   {u \cdot {\prod\limits_{i =
1}^{n} {{\left[ {x_{i}}  \right]}\,}} r^{ - 2\tilde {\alpha} _{n}
}F_{A}^{\left( {n} \right)} {\left[ {{\begin{array}{*{20}c}
 {\tilde {\alpha} _{n} ,1 - \alpha _{1} ,...,1 - \alpha _{n} ;} \hfill \\
 {2 - 2\alpha _{1} ,...,2 - 2\alpha _{n} ;} \hfill \\
\end{array}} \sigma}  \right]} \cdot {\sum\limits_{i = 1}^{n} {{\frac{{1 -
2\alpha _{i}}} {{x_{i}}} }\cos \left( {{\rm {\bf n}};x_{i}}  \right)}
}dC_{\rho}} }  .
\]

We use the following generalization spherical system of coordinates:

\[
x_{1} = \xi _{1} + \rho \cos \varphi _{1} ,
\]

\[
x_{2} = \xi _{2} + \rho \sin \varphi _{1} \cos \varphi _{2} ,
\]

\[
x_{3} = \xi _{3} + \rho \sin \varphi _{1} \sin \varphi _{2} \cos \varphi _{3} ,
\]

\[
................................
\]

\[
x_{m - 1} = \xi _{m - 1} + \rho \sin \varphi _{1} \sin \varphi _{2} ...\sin \varphi
_{m - 2} \cos \varphi _{m - 1} ,
\]

\[
x_{m} = \xi _{m} + \rho \sin \varphi _{1} \sin \varphi _{2} ...\sin \varphi _{m - 2}
\sin \varphi _{m - 1}
\]

\[
{\left[ {\rho \ge 0,\,\,0 \le \varphi _{1} \le \pi ,\,...,\,0 \le \varphi _{m - 2}
\le \pi ,\,\,\,0 \le \varphi _{m - 1} \le 2\pi \,\,} \right]}.
\]

Then we have

\[
I_{11} = 2\tilde {\alpha} _{n} \gamma _{n} {\prod\limits_{i = 1}^{n}
{{\left[ {\xi _{i}^{1 - 2\alpha _{i}}}   \right]}}} \,\rho ^{2\alpha _{1} +
... + 2\alpha _{n} - 2n}{\int\limits_{0}^{2\pi}  {d\varphi _{m - 1}}
}{\int\limits_{0}^{\pi}  {\sin \varphi _{m - 2} d\varphi _{m - 2}}
}...{\int\limits_{0}^{\pi}  {\tilde {u}\,v\,F_{A}^{\left( {n} \right)}
\left( {\tilde {\sigma}}  \right)\sin ^{m - 2}\phi _{1} d\varphi _{1}}}  ,
\]

\noindent
where

\[
\begin{array}{l}
 \tilde {u}: = u\left( {\xi _{1} + \rho \cos \varphi _{1} ,\xi _{2} + \rho \sin
\varphi _{1} \cos \varphi _{2} ,\xi _{3} + \rho \sin \varphi _{1} \sin \varphi _{2} \cos
\varphi _{3} ,...,} \right. \\
\,\,\,\,\left. {\xi _{m - 1} + \rho \sin \varphi _{1} \sin \varphi _{2} ...\sin
\varphi _{m - 2} \cos \varphi _{m - 1} ,\xi _{m} + \rho \sin \varphi _{1} \sin \varphi
_{2} ...\sin \varphi _{m - 2} \sin \varphi _{m - 1}}  \right), \\
 \end{array}
\]

\[
\begin{array}{l}
 v: = \left( {\xi _{1} + \rho \cos \varphi _{1}}  \right)\left( {\xi _{2} +
\rho \sin \varphi _{1} \cos \varphi _{2}}  \right)\left( {\xi _{3} + \rho \sin
\varphi _{1} \sin \varphi _{2} \cos \varphi _{3}}  \right)... \\
\,\,\,\,...\left( {\xi _{m - 1} + \rho \sin \varphi _{1} \sin \varphi _{2}
...\sin \varphi _{m - 2} \cos \varphi _{m - 1}}  \right)\left( {\xi _{m} + \rho
\sin \varphi _{1} \sin \varphi _{2} ...\sin \varphi _{m - 2} \sin \varphi _{m - 1}}
\right), \\
 \end{array}
\]

\[
F_{A}^{\left( {n} \right)} \left( {\tilde {\sigma}}   \right): =
F_{A}^{\left( {n} \right)} \left( {\tilde {\alpha} _{n} + 1,1 - \alpha _{1}
,...,1 - \alpha _{n} ;2 - 2\alpha _{1} ,...,2 - 2\alpha _{n} ;\tilde {\sigma
}}   \right),
\]

\[
{\rm \tilde {\sigma}}  : = \left( {\sigma _{1\rho}  ,...,\sigma
_{n\rho}}   \right);
r_{k}^{2} = \left( {x_{k} + \xi _{k}}  \right)^{2} + {\sum\limits_{i = 1,i
\ne k}^{m} {\left( {x_{i} - \xi _{i}}  \right)^{2}}} ,
\quad
\sigma _{k\rho}  = 1 - {\frac{{r_{k\rho} ^{2}}} {{\rho ^{2}}}}.
\]

First we evaluate $F_{A}^{\left( {n} \right)} \left( {\tilde {\sigma}}  \right)$. For this aim we use decomposition formula (\ref{eq12}) and then auto-transformation formula (\ref{eq6}):

\[
F_{A}^{\left( {n} \right)} \left( {\tilde {\sigma}}   \right) =
{\sum\limits_{{\mathop {m_{i,j} = 0}\limits_{(2 \le i \le j \le n)}
}}^{\infty}  {{\frac{{(\tilde {\alpha} _{n} + 1)_{N(n,n)}}} {{m_{ij} !}}}}
}{\prod\limits_{k = 1}^{n} {{\left[ {{\frac{{(1 - \alpha _{k} )_{M(k,n)}
}}{{(2 - 2\alpha _{k} )_{M(k,n)}}} }\left( {{\frac{{r_{k\rho} ^{2}}} {{\rho
^{2}}}}} \right)^{\alpha _{k} - 1 - M(k,n)}\left( {1 - {\frac{{r_{k\rho
}^{2}}} {{\rho ^{2}}}}} \right)^{M(k,n)}} \right]}}}
\]

\[
\times {\prod\limits_{k = 1}^{n} {{\left[ {F\left( {1 - 2\alpha _{k} -
\tilde {\alpha} _{n} + M(k,n) - N(k,n),1 - \alpha _{k} + M(k,n);2 - 2\alpha
_{k} + M(k,n);{\frac{{\sigma _{k\rho}} } {{\sigma _{k\rho}  - 1}}}} \right)}
\right]}}} ,
\]

\noindent
where $M(k,n)$ and $N(k,n)$ are an expressions defined in (\ref{eq13}).

After the elementary evaluations we find

\begin{equation}
\label{eq29}
F_{A}^{\left( {n} \right)} \left( {\tilde {\sigma}}   \right) = \rho
^{2n - 2\alpha _{1} - ... - 2\alpha _{n}} {\prod\limits_{k = 1}^{n} {{\left[
{r_{k\rho} ^{2\alpha _{k} - 2}}  \right]}}} \, \cdot \aleph ,
\end{equation}

\noindent
where

\[
\aleph : = {\sum\limits_{{\mathop {m_{i,j} = 0}\limits_{(2 \le i \le j \le
n)}}} ^{\infty}  {{\frac{{(\tilde {\alpha} _{n} + 1)_{N(n,n)}}} {{m_{ij}
!}}}{\prod\limits_{k = 1}^{n} {{\left[ {{\frac{{(1 - \alpha _{k} )_{M(k,n)}
}}{{(2 - 2\alpha _{k} )_{M(k,n)}}} }\left( {{\frac{{\rho ^{2}}}{{r_{k\rho
}^{2}}} } - 1} \right)^{M(k,n)}} \right]}}}} }
\]

\[
\times {\prod\limits_{k = 1}^{n} {{\left[ {F\left( {1 - 2\alpha _{k} -
\tilde {\alpha} _{n} + M(k,n) - N(k,n),1 - \alpha _{k} + M(k,n);2 - 2\alpha
_{k} + M(k,n);1 - {\frac{{\rho ^{2}}}{{r_{k\rho} ^{2}}} }} \right)}
\right]}}} .
\]

It is easy to see that when $\rho \to 0$ the function $\aleph $ becomes an
expression that does not depend on ${\rm x}$ and ${\rm \xi} $. Indeed,
taking into account the equality (\ref{eq14}), we have

\[
{\mathop {\lim} \limits_{\rho \to 0}} \aleph = {\sum\limits_{{\mathop
{m_{i,j} = 0}\limits_{(2 \le i \le j \le n)}}} ^{\infty}  {{\frac{{(\tilde
{\alpha} _{n} + 1)_{N(n,n)}}} {{m_{ij} !}}}}} {\prod\limits_{k = 1}^{n}
{{\left[ {{\frac{{(1 - \alpha _{k} )_{M(k,n)}}} {{(2 - 2\alpha _{k}
)_{M(k,n)}}} }} \right]}}}
\]

\begin{equation}
\label{eq30}
\times {\prod\limits_{k = 1}^{n} {{\left[ {F\left( {1 - 2\alpha _{k} -
\tilde {\alpha} _{n} + M(k,n) - N(k,n),1 - \alpha _{k} + M(k,n);2 - 2\alpha
_{k} + M(k,n);1} \right)} \right]}}} .
\end{equation}

Applying now the summation formula (\ref{eq5}) to each hypergeometric function
$F\left( {a,b;c;1} \right)$ in the sum (\ref{eq30}), we get

\[
{\mathop {\lim} \limits_{\rho \to 0}} \aleph = {\frac{{1}}{{\Gamma \left(
{\tilde {\alpha} _{n} + 1} \right)}}}{\sum\limits_{{\mathop {m_{i,j} =
0}\limits_{(2 \le i \le j \le n)}}} ^{\infty}  {{\frac{{\Gamma \left(
{\tilde {\alpha} _{n} + 1 + N(n,n)} \right)}}{{m_{ij} !}}}}
}
\]
\[
{\times\prod\limits_{k = 1}^{n} {{\left[ {{\frac{{\Gamma (1 - \alpha _{k} +
M(k,n))\Gamma \left( {2 - 2\alpha _{k}}  \right)\Gamma \left( {\tilde
{\alpha} _{n} + \alpha _{k} + N(k,n) - M(k,n)} \right)}}{{\Gamma \left(
{\tilde {\alpha} _{n} + 1 + N(k,n)} \right)\Gamma ^{2}\left( {1 - \alpha
_{k}}  \right)}}}} \right]}}} .
\]

Taking into account the identity

\begin{equation}
\label{eq31}
{\sum\limits_{{\mathop {m_{i,j} = 0}\limits_{(2 \le i \le j \le n)}
}}^{\infty}  {{\frac{{\left( {\tilde {\alpha} _{n} + 1} \right)_{N(n,n)}
}}{{m_{ij} !}}}}} {\prod\limits_{k = 1}^{n} {{\left[ {{\frac{{\left( {\tilde
{\alpha} _{n} + \alpha _{k}}  \right)_{N(k,n) - M(k,n)} (1 - \alpha _{k}
)_{M(k,n)}}} {{\left( {\tilde {\alpha} _{n} + 1} \right)_{N(k,n)}}} }}
\right]}}}  = \Gamma \left( {{\frac{{m}}{{2}}}} \right){\prod\limits_{k =
1}^{n} {\Gamma ^{ - 1}\left( {\tilde {\alpha} _{n} + \alpha _{k}}  \right)}
},
\end{equation}

\noindent
we obtain

\begin{equation}
\label{eq32}
{\mathop {\lim} \limits_{\rho \to 0}} \aleph = {\frac{{\Gamma \left( {m / 2}
\right)}}{{\Gamma \left( {\tilde {\alpha} _{n} + 1}
\right)}}}{\prod\limits_{i = 1}^{n} {{\frac{{\Gamma \left( {2 - 2\alpha _{k}
} \right)}}{{\Gamma \left( {1 - \alpha _{k}}  \right)}}}}} .
\end{equation}

Using the properties (\ref{eq3}) of gamma-function $\Gamma \left( {z} \right)$,
property (\ref{eq4}) of the Pochhammer symbol and summation formula (\ref{eq5}) for
hypergeometric function $F\left( {a,b;c;z} \right)$, the formula (\ref{eq31}) is
proved by the method mathematical induction.

Now we consider an integral

\[
L_{m} = {\int\limits_{0}^{2\pi}  {d\varphi _{m - 1}}}  {\int\limits_{0}^{\pi}
{\sin \varphi _{m - 2} d\varphi _{m - 2}}}  {\int\limits_{0}^{\pi}  {\sin ^{2}\varphi
_{m - 3} d\varphi _{m - 3}}}  ...{\int\limits_{0}^{\pi}  {\sin ^{m - 2}\varphi
_{1} d\varphi _{1}}}  .
\]

\noindent
with elementary transformations it is not difficult to establish that

\begin{equation}
\label{eq33}
L_{2m} = {\frac{{2\,\pi ^{m}}}{{(m - 1)!}}},\,\,L_{2m + 1} = {\frac{{2^{m +
1}\,\pi ^{m}}}{{(2m - 1)!!}}},\,\,\,m = 1,2,3,...
\end{equation}

If we take into account (\ref{eq29}), (\ref{eq32}), (\ref{eq33}) and (\ref{eq16}), then we will have

\[
{\mathop {\lim} \limits_{\rho \to 0}} I_{11} = u\left( {\xi}  \right).
\]

By similar evaluations one can get that

\[
{\mathop {\lim} \limits_{\rho \to 0}} I_{12} = {\mathop {\lim} \limits_{\rho
\to 0}} I_{13} = {\mathop {\lim} \limits_{\rho \to 0}} I_{2} = 0.
\]

If we consider an integral

\[
{\int_{C_{\rho}}   {{\rm x}^{\left( {2\alpha}  \right)}G\left( {{\rm x};{\rm
\xi}}  \right){\frac{{\partial u\left( {{\rm x}} \right)}}{{\partial {\rm
{\bf n}}}}}dC_{\rho}} }  ,
\]

\noindent
using above given algorithm for evaluations (in this case calculations will
be more simple), we can prove that

\[
{\mathop {\lim} \limits_{\rho \to 0}} {\int_{C_{\rho}}   {{\rm x}^{\left(
{2\alpha}  \right)}G\left( {{\rm x};{\rm \xi}}  \right){\frac{{\partial
u\left( {{\rm x}} \right)}}{{\partial {\rm {\bf n}}}}}dC_{\rho}} }   = 0.
\]

Now from (\ref{eq23}) we can write the solution of the Dirichlet problem as
follows:

\begin{equation}
\label{eq34}
u\left( {\xi}  \right) = {\sum\limits_{k = 1}^{n} {{\int_{S_{k}}  {\tilde
{x}_{k}^{\left( {2\alpha}  \right)} G_{k}^{\ast}  \left( {x_{1} ,...,x_{k -
1} ,0,x_{k + 1} ,...,x_{m} ;{\rm \xi}}  \right)\tau _{k} \left( {\tilde
{x}_{k}}  \right)dS_{k}}}  \,}}  + {\int_{S} {x^{\left( {2\alpha}
\right)}{\frac{{\partial G\left( {{\rm x};{\rm \xi}}
\right)}}{{\partial {\rm {\bf n}}}}}\phi \left( {{\rm x}} \right)dS}} .
\end{equation}

The particular values of Green's function are given by

\[
G_{k}^{\ast}  \left( {x_{1} ,...,x_{k - 1} ,0,x_{k + 1} ,...,x_{m} ;{\rm \xi
}} \right) = \left( {1 - 2\alpha _{k}}  \right)\gamma _{n} \xi _{k}^{1 -
2\alpha _{k}}  {\prod\limits_{i = 1,i \ne k}^{n} {{\left[ {\left( {x_{i} \xi
_{i}}  \right)^{1 - 2\alpha _{i}}}  \right]}\,}}
\]

\[
\times {\left\{ {{\frac{{F_{A}^{\left( {n - 1} \right)} {\left[
{{\begin{array}{*{20}c}
 {\tilde {\alpha} _{n} ,1 - \alpha _{1} ,...,1 - \alpha _{k - 1} ,1 - \alpha
_{k + 1} ,...,1 - \alpha _{n} ;} \hfill \\
 {2 - 2\alpha _{1} ,...,2 - 2\alpha _{k - 1} ,2 - 2\alpha _{k + 1} ,...,2 -
2\alpha _{n} ;} \hfill \\
\end{array}} \sigma _{0}}  \right]}}}{{{\left[ {\xi _{k}^{2} +
{\sum\limits_{i = 1,i \ne k}^{m} {\left( {\xi _{i} - x_{i}}  \right)^{2}}}}
\right]}^{\,\tilde {\alpha} _{n}}} }}} \right.}
\]

\[
{\left. { - \left( {{\frac{{a}}{{R_{0}}} }} \right)^{2n - 4\alpha _{1} - ...
- 4\alpha _{n}} {\frac{{F_{A}^{\left( {n - 1} \right)} {\left[
{{\begin{array}{*{20}c}
 {\tilde {\alpha} _{n} ,1 - \alpha _{1} ,...,1 - \alpha _{k - 1} ,1 - \alpha
_{k + 1} ,...,1 - \alpha _{n} ;} \hfill \\
 {2 - 2\alpha _{1} ,...,2 - 2\alpha _{k - 1} ,2 - 2\alpha _{k + 1} ,...,2 -
2\alpha _{n} ;} \hfill \\
\end{array}} \bar {\sigma} _{0}}  \right]}}}{{{\left[ {{\sum\limits_{i = 1,i
\ne k}^{m} {\left( {a - {\frac{{x_{i} \xi _{i}}} {{a}}}} \right)^{2}}}  +
{\frac{{1}}{{a^{2}}}}{\sum\limits_{i = 1,i \ne k}^{m} {{\sum\limits_{j = 1,j
\ne i}^{m} {x_{i}^{2} \xi _{j}^{2}}} } }  - \left( {m - 2} \right)a^{2}}
\right]}^{\,\tilde {\alpha} _{n}}} }}} \right\}} ,
\]

\bigskip

\noindent
where

\[
\sigma _{0} : = \left( {\sigma _{1}^{0} ,...,\sigma _{k - 1}^{0} ,\sigma _{k
+ 1}^{0} ,...,\sigma _{n}^{0}}  \right),
\quad
\bar {\sigma} _{0} : = \left( {\bar {\sigma} _{1}^{0} ,...,\bar {\sigma} _{k
- 1}^{0} ,\bar {\sigma} _{k + 1}^{0} ,...,\bar {\sigma} _{n}^{0}}  \right),
\]

\[
\sigma _{s}^{0} = - {\frac{{4x_{s} \xi _{s}}} {{\xi _{k}^{2} +
{\sum\limits_{i = 1,i \ne k}^{m} {\left( {\xi _{i} - x_{i}}  \right)^{2}}
}}}},
\]

\[
\bar {\sigma} _{s}^{0} = - {\frac{{a^{2}}}{{R_{0}^{2}}} }{\frac{{4x_{s} \xi
_{s}}} {{{\sum\limits_{i = 1,i \ne k}^{m} {\left( {a - {\frac{{x_{i} \xi
_{i}}} {{a}}}} \right)^{2}}}  + {\frac{{1}}{{a^{2}}}}{\sum\limits_{i = 1,i
\ne k}^{m} {{\sum\limits_{j = 1,j \ne i}^{m} {x_{i}^{2} \xi _{j}^{2}}} } }
- \left( {m - 2} \right)a^{2}}}},
\quad
s = \overline {1,n} ,\,\,s \ne k.
\]

Constant $\gamma _{n} $ has the form (\ref{eq16}).

Hence, the main result of the paper is formulated as the following theorem:

\textbf{Theorem.} If $\tau _{k} \left( {\tilde {x}_{k}}  \right) \in
C^{2}\left( {S_{k}}  \right)$ and $\varphi \left( {{\rm x}} \right) \in
C^{2}\left( {S} \right)$ are given functions fulfilling the matching
conditions (\ref{eq19}), then the Dirichlet problem has unique solution represented
by formula (\ref{eq34}).

\bigskip

Ergashev Tuhtasin Gulamjanovich

Doctoral student, Uzbekistan Academy of Sciences, V.I.Romanovskiy Institute
of Mathematics,

81, Mirzo Ulugbek street, Tashkent, 100170 \\
ergashev.tukhtasin@gmail.com

\end{document}